\documentclass[12pt]{amsart}
\usepackage{amscd,amssymb}
\usepackage[graph,frame,poly,arc]{xy}  
\usepackage[plainpages,backref,urlcolor=blue]{hyperref}

\theoremstyle{plain}
\newtheorem{thm}[subsection]{Theorem}
\newtheorem{lem}[subsection]{Lemma}
\newtheorem{prop}[subsection]{Proposition}
\newtheorem{cor}[subsection]{Corollary}

\theoremstyle{definition}
\newtheorem{rk}[subsection]{Remark}
\newtheorem{definition}[subsection]{Definition}
\newtheorem{ex}[subsection]{Example}
\newtheorem{conj}[subsection]{Conjecture}

\numberwithin{equation}{section}
\newcommand{\PPP}{{\mathcal P}}

\newcommand{\A}{{\mathcal A}}

\newcommand{\al}{{\alpha}}
\newcommand{\be}{{\beta}}

\newcommand{\Z}{\mathbb{Z}}
\newcommand{\Q}{\mathbb{Q}}

\newcommand{\C}{\mathbb{C}}
\newcommand{\PP}{\mathbb{P}}

\DeclareMathOperator{\dd}{d}
\begin{document}
%\date{June 4, 2009}

\title [A computational approach to Milnor fiber cohomology]
{A computational approach to Milnor fiber cohomology }

\author[Alexandru Dimca]{Alexandru Dimca$^1$}
\address{Univ. Nice Sophia Antipolis, CNRS,  LJAD, UMR 7351, 06100 Nice, France. }
\email{dimca@unice.fr}

\author[Gabriel Sticlaru]{Gabriel Sticlaru}
\address{Faculty of Mathematics and Informatics,
Ovidius University,
Bd. Mamaia 124, 
900527 Constanta,
Romania}
\email{gabrielsticlaru@yahoo.com }
\thanks{$^1$ Partially supported by Institut Universitaire de France.}

\subjclass[2010]{Primary 32S55; Secondary 32S35, 32S22.}

\keywords{plane curve, Milnor fiber, monodromy}

\begin{abstract} In this note we consider the Milnor fiber $F$ associated to a reduced projective plane curve $C$. A computational approach for the determination of the characteristic polynomial of the monodromy action on the first cohomology group of $F$, also known as the Alexander polynomial of the curve $C$, is presented. This leads to an effective algorithm to detect all the monodromy eigenvalues and, in many cases, explicit bases for the monodromy eigenspaces in terms of polynomial differential forms.
\end{abstract}
 
\maketitle

%\tableofcontents

\section{Introduction} \label{sec:intro}

Let $C:f=0$ be a reduced plane curve  of degree $d\geq 3$  in the complex projective plane $\PP^{2}$, defined by a homogeneous polynomial $f \in S=\C[x,y,z]$.
Consider the corresponding complement $U=\PP^{2}\setminus C$, and the global Milnor fiber $F$ defined by $f(x,y,z)=1$ in $\C^3$ with monodromy action $h:F \to F$, $h(x)=\exp(2\pi i/d)\cdot (x,y,z)$. 
One can consider  the characteristic polynomials of the monodromy, namely
\begin{equation} 
\label{Delta}
\Delta^j_C(t)=\det (t\cdot Id -h^j|H^j(F,\C)),
\end{equation} 
for $j=0,1,2$. It is clear that, when the curve $C$ is reduced,  one has $\Delta^0_C(t)=t-1$, and moreover
\begin{equation} 
\label{Euler}
\Delta^0_C(t)\Delta^1_C(t)^{-1}\Delta^2_C(t)=(t^d-1)^{\chi(U)},
\end{equation} 
where $\chi(U)$ denotes the Euler characteristic of the complement $U$, see for instance \cite[Proposition 4.1.21]{D1}.
It follows that the polynomial $\Delta_C(t)=\Delta^1_C(t)$, also called the Alexander polynomial of $C$, see \cite{R}, determines the remaining polynomial $\Delta^2_C(t)$. To determine the Alexander polynomial $\Delta_C(t)$, or equivalently, the eigenvalues of the monodromy operator
\begin{equation} 
\label{mono1}
h^1: H^1(F,\C) \to H^1(F,\C)
\end{equation} 
starting from $C$ or  $f$  is a rather difficult problem, going back to O. Zariski and attracting an extensive literature, see for instance \cite{HE}, \cite{L1}, \cite{L2}, \cite{OkaS}, \cite{K},\cite{Deg}, \cite{AD},  \cite{D1}.

Let $\Omega^j$ denote the graded $S$-module of (polynomial) differential $j$-forms on $\C^3$, for $0 \leq j \leq 3$. The complex $K^*_f=(\Omega^*, \dd f \wedge)$ is nothing else but the Koszul complex in $S$ of the partial derivatives $f_x$, $f_y$ and $f_z$ of the polynomial $f$. The general theory says that there is a spectral sequence $E_*(f)$, whose first term $E_1(f)$ is computable from the cohomology of the Koszul complex $K^*_f$ and whose limit
$E_{\infty}(f)$ gives us the action of monodromy operator  on the graded pieces of $H^*(F,\C)$ with respect to a filtration $P$ to be defined below, see \cite{Dcomp}, \cite[Chapter 6]{D1}, \cite{DS1}.
Moreover, note that the $1$-eigenspace $H^1(F,\C)_1$ coincides with $H^1(U,\C)$, and an explicit basis for this cohomology group can be easily given using differential $1$-forms, see \cite[Example (4.1)]{Dcomp}. The sum of the other eigenspaces, denoted by $H^1(F,\C)_{\ne1}$,
is a pure Hodge structure of weight $1$, see \cite{DP}, and the study of the spectral sequence
 $E_*(f)$ often provides a basis for $H^1(F,\C)_{\ne1}$, given in terms of differential $1$-forms and adapted to the Hodge structure.

This approach is particularly useful when $C$ is assumed to have only weighted homogeneous singularities, e.g.  when $C$ is a line arrangement $\A$,
in view of the following conjecture, see for a more precise statement Conjecture  \ref{conj} below.

\begin{conj}
\label{conj0}
If the reduced plane curve $C:f=0$ has only weighted homogeneous singularities, then
the spectral sequence $E_*(f)$ degenerates at the $E_2$-term. 
\end{conj}
The converse implication is known to hold, see \cite[Theorem 5.2]{DS1}.
The first aim of this note is to test  Conjecture  \ref{conj0} in a large number of situations and to show how  the spectral sequence $E_*(f)$ can be applied to obtained explicit bases for the eigenspaces of  monodromy operators \eqref{mono1} of line arrangements, see Theorem \ref{thm1} and  Theorem \ref{thmHA}, or for some irreducible plane curves, see Propositions \ref{propNWH}, \ref{prop9cusp}, \ref{propZariski6}, \ref{propZariski12}

When $G$ is a finite group acting on the Milnor fiber $F$, results as those listed above allows us to determine the action of $G$ on the cohomology group $H^1(F,\C)$, see Examples \ref{exaction1}, 
\ref{exaction2}.

The second aim of this note is to show that, even when the spectral sequence $E_*(f)$ does not degenerate at $E_2$, a small set of its second order terms $E_2^{s,t}(f)$ can be used to detect all the eigenvalues of the monodromy operator \eqref{mono1}, see next section for the notation used in the following, which is the main result of our note.

\begin{thm}
\label{thmEV}
Let $C:f=0$ be a reduced degree $d$ curve in $\PP^2$, and let $$\lambda=\exp (-2\pi i k/d)\ne 1,$$  with $k \in (0,d)$ an integer. 
Then $\lambda$ is a root of the Alexander polynomial $\Delta_C(t)$ if and only if
either
$E_2^{1,0}(f)_k \ne 0$ or $E_2^{1,0}(f)_{k'} \ne 0$, where $k'=d-k$. 
The multiplicity $m(\lambda)$
of the root $\lambda$ satisfies the inequalities
$$\max (\dim E_2^{1,0}(f)_{k}, \dim E_2^{1,0}(f)_{k'}) \leq m(\lambda) \leq \dim E_2^{1,0}(f)_{k}+ \dim E_2^{1,0}(f)_{k'}.   $$
In particular, equality holds everywhere if either $ E_2^{1,0}(f)_{k}=0$ or  $E_2^{1,0}(f)_{k'}=0$.
\end{thm}
The Hessian line arrangement discussed in Theorem \ref{thmHA} shows that the first inequality in Theorem \ref{thmEV} can be strict. On the other hand, for all the examples discussed in this note, the second inequality is in fact an equality.

As an application, consider the family of plane curves
\begin{equation} 
\label{curve}
C_{5m}: f_{5m}=(y^mz^m-x^{2m})^2y^m-x^{5m},
\end{equation} 
for $m \geq 1$ introduced in \cite{B+}. These curves are free with exponents $(2m,3m-1)$, and for $m$ odd they give examples of free irreducible curves which are not rational, see \cite[Theorem 3.9]{B+}. When $m$ is even, these curves have two irreducible components and a large Alexander polynomial as shown by the following result, which is proved in subsection \ref{freecurve}.

\begin{cor}
\label{corAP}
For $2 \leq m \leq 20$, $m$  even, the Alexander polynomial of the curve $C_{5m}$ defined in  \eqref{curve} is given by
$$\Delta_{C_{5m}}(t)=t^{5m}-1.$$
On the other hand, for  $1 \leq m \leq 19$, $m$  odd, the Alexander polynomial of the curve $C_{5m}$ defined in \eqref{curve} is trivial, i.e. $\Delta_{C_{5m}}(t)=1.$
\end{cor}
We expect the above result to hold for any $m$, but such a claim is beyond our computational approach. For $m \geq 2$, the curve $C_{5m}$ has some non weighted homogeneous singularities and hence the spectral sequence $E_*(f)$ does not degenerate at $E_2$.

The computations in this note were made using two computer algebra systems, namely CoCoA \cite{Co} and Singular \cite{Sing}.
The corresponding codes are available on request.

\section{Gauss-Manin complexes, Koszul complexes, and Milnor fiber cohomology} \label{sec2}

Let $S$ be the polynomial ring $\C[x,y,z]$ with the usual grading and consider a reduced homogeneous  polynomial $f \in S$ of degree $d$. The graded Gauss-Manin complex $C_f^*$ associated to $f$ is defined by taking $C_f^j=\Omega^j[\partial_t]$, i.e. formal polynomials in $\partial_t$ with coefficients in the space of differential forms $\Omega^j$, where $\deg \partial_t=-d$ and the differential
$\dd: C_f^j \to C_f^{j+1}$ is $\C$-linear and given by 
\begin{equation} 
\label{difC}
 \dd (\omega \partial_t^q)=(\dd \omega)\partial_t^q-(\dd f \wedge \omega) \partial_t^{q+1},
\end{equation} 
see for more details \cite[Section 4]{DS1}. 
The complex $C_f^*$ has a natural increasing filtration $P'_*$ defined by
\begin{equation} 
\label{filC}
 P'_qC^j_f=\oplus_{i \leq q+j}\Omega^j\partial_t^i.
\end{equation} 
If we set $P'^q=P'_{-q}$ in order to get a decreasing filtration, then one has
\begin{equation} 
\label{grC}
 Gr^q_{P'}C^*_f=\sigma_{\geq q}(K^*_f((3-q)d)),
\end{equation} 
the truncation of a shifted version of the Koszul complex $K^*_f$.
Moreover, this yields a decreasing filtration $P'$ on the cohomology groups $H^j(C^*_f)$ and a spectral sequence
\begin{equation} 
\label{spsqC}
 E_1^{q,j-q}(f) \Rightarrow H^j(C^*_f).
\end{equation} 
On the other hand, the cohomology $H^j(F,\C)$ of the Milnor fiber $F:f(x,y,z)=1$ associated to $f$ has a pole order decreasing filtration $P$, see \cite[Section 3]{DS1}, such that there is a natural identification for any integers $q$, $j$ and $k \in [1,d]$
\begin{equation} 
\label{filH}
 P'^{q+1}H^{j+1}(C^*_f)_k=P^qH^j(F,\C)_{\lambda},
\end{equation} 
where $\lambda=\exp (-2 \pi ik/d).$ Moreover, the $E_1$-term of the spectral sequence 
\eqref{spsqC} is completely determined by the morphism of graded $\C$-vector spaces
\begin{equation} 
\label{differential1}
 \dd ' : H^2(K^*_f) \to H^3(K^*_f),
\end{equation} 
induced by the exterior differentiation of forms, i.e. $\dd ' :[\omega] \mapsto [\dd (\omega)]$.
Note that this morphism $\dd'$ coincides with the morphism $\dd ^{(1)}: N \to M$ considered in 
\cite{DS1}, up-to some shifts in gradings. More precisely, if we set for $k \in [1,d]$,
\begin{equation} 
\label{newspsq}
E_1^{s,t}(f)_k=H^{s+t+1}(K^*_f)_{td+k},
\end{equation} 
we get a new form of (a homogeneous component of) the spectral sequence \eqref{spsqC}, and hence with the above notation one has
\begin{equation} 
\label{limit}
E_{\infty}^{s,t}(f)_k=Gr_P^sH^{s+t}(F,\C)_{\lambda }.
\end{equation} 
The case $k=d$ is also discussed in \cite{DStEdin}.
We need the following result, see  \cite[Theorem 5.3]{DS1}.

\begin{thm}
\label{thminj}
Assume that the plane curve $C:f=0$ in $\PP^2$ has only weighted homogeneous singularities.
Then the dimension of 
$$E_2^{1-t,t}(f)_k=\ker \{ \dd ' : H^2(K^*_f)_{td+k} \to H^3(K^*_f)_{td+k}\}$$
is upper bounded by the number of spectral numbers $\al_{p_i,j}$ of the singularity $(C,p_i)$ equal to $\frac{td+k}{d}$,
when $p_i$ ranges over all the singularities of $C$ and $j=1,...,\mu(C,p_i).$
\end{thm}

\begin{ex}
\label{exspec} For the monomial arrangement $\A(m,m,3)$ defined by 
$$\A(m,m,3): f=(x^m-y^m)(x^m-z^m)(y^m-z^m)=0$$
for $m\geq 2$, there are two types of singularities, the triple points, with local equation $u^3+v^3=0$ and hence with maximal spectral number 
$\al_3=4/3$ and the three $m$-fold intersection points, e.g. $[1:0:0]$, with local equation $u^m+v^m=0$ and maximal spectral number 
$\al_m=\frac{2m-2}{m}$, see \cite{DS14} for such computations of spectral numbers. In the case $m \geq 3$, it follows that 
$$ \dd '_l : H^2(K^*_f)_l \to H^3(K^*_f)_l$$
 is injective for $l >6m-6=2(d-3)$.

\end{ex}

One has the following, see Question 2 in \cite{DS1} for a more general setting.

\begin{conj}
\label{conj}
The spectral sequences \eqref{spsqC} and \eqref{newspsq} degenerate at the $E_2$-term when the reduced plane curve $C:f=0$ has only weighted homogeneous singularities.
\end{conj}

\begin{rk}
\label{rkNWH}
Even in the case when the spectral sequence \eqref{newspsq} does not degenerate at the $E_2$-term we have 
\begin{equation} 
\label{degNWH}
E_2^{1-t,t}(f)_k=E_{\infty}^{1-t,t}(f)_k=Gr_P^{1-t}H^{1}(F,\C)_{\lambda },
\end{equation} 
for $t=0$, or for $t=1$ and $k<3$, since the second differential $$d_2: E_2^{1-t,t}(f)_k \to E_2^{3-t,t-1}(f)_k$$
vanishes.
Indeed, one has
\begin{equation} 
\label{deg2}
E_1^{3-t,t-1}(f)_k=H^3(K_f^*)_{(t-1)d+k}=0
\end{equation} 
for $(t-1)d+k<3$, i.e. for $t=0$ and any $k \in [1,d]$, or for $t=1$ and $k<3$. 
\end{rk}

\subsection{Proof of Theorem \ref{thmEV}}

By \eqref{degNWH}, we have the following
$$E_2^{1,0}(f)_k=E_{\infty}^{1,0}(f)_k=Gr_P^{1}H^{1}(F,\C)_{\lambda },$$
and
$$E_2^{1,0}(f)_{k'}=E_{\infty}^{1,0}(f)_{k'}=Gr_P^{1}H^{1}(F,\C)_{\overline \lambda }.$$
If $E_2^{1,0}(f)_k \ne 0$, it follows that $\lambda$ is an eigenvalue of the monodromy operator.
If $E_2^{1,0}(f)_{k'} \ne 0$, it follows that $\overline \lambda$ is an eigenvalue of the monodromy operator. But this operator is defined over $\Q$, hence 
$m(\lambda)=m(\overline\lambda) $. 

Note that on the cohomology group $H^1(F,\C)$ one has $P^2=0$ and $F^1 \subset P^1$, where $F^1$ denotes the Hodge filtration. Let now $\lambda$ be an eigenvalue of the monodromy operator.
The fact that $H^1(F,\C)_{\ne 1}$ is a pure Hodge structure of weight 1, implies that either
$H^{1,0}(F,\C)_{\lambda } \ne 0$ or $H^{0,1}(F,\C)_{\lambda } \ne 0$. In the first case we get
$Gr_P^{1}H^{1}(F,\C)_{\lambda }\ne 0$, while in the second case we get $H^{1,0}(F,\C)_{\overline \lambda } \ne 0$ and hence $Gr_P^{1}H^{1}(F,\C)_{\overline \lambda }\ne 0$.
This proves the first claim as well as the first inequality in 
Theorem \ref{thmEV}. 

\medskip

For the second inequality, note that one has
$$\dim E_2^{1,0}(f)_{k}+ \dim E_2^{1,0}(f)_{k'} =\dim P^{1}H^{1}(F,\C)_{\lambda }+\dim P^{1}H^{1}(F,\C)_{\overline \lambda }\geq $$
$$\geq \dim F^{1}H^{1}(F,\C)_{\lambda }+\dim F^{1}H^{1}(F,\C)_{\overline \lambda }=\dim H^{1,0}(F,\C)_{\lambda }+\dim H^{1,0}(F,\C)_{\overline \lambda }=$$
$$=\dim H^{1,0}(F,\C)_{\lambda }+\dim H^{0,1}(F,\C)_{ \lambda }=\dim H^{1}(F,\C)_{\lambda } = m(\lambda).$$

This completes the proof of Theorem \ref{thmEV}.

\begin{rk}
\label{rkcomplex}

(i) It is sometimes useful to replace the Gauss-Manin complex $C^*_f$ by the simpler complex $(\Omega^*,D_f)$, where the differential $D_f$ is given by
\begin{equation} 
\label{newcomplex}
D_f(\omega)=\dd \omega -\frac{|\omega|}{d} \dd f \wedge \omega,
\end{equation} 
with $\omega \in \Omega ^j$ a homogeneous differential forms of degree $|\omega|$, see \cite{D1}, p. 190. This complex has the advantage that one has explicit isomorphisms 
\begin{equation} 
\label{iso1}
H^j(\Omega^*,D_f)= H^{j-1}(F,\C),  \  \  \    [\omega]\mapsto \iota^*(\Delta( \omega)),
\end{equation} 
for any $j$, where $\iota: F \to \C^3$ denotes the inclusion of the Milnor fiber $F$ into $\C^3$ and $\Delta: \Omega^j \to \Omega ^{j-1}$ denotes the contraction with the Euler vector field , see \cite{D1}, p. 193.

\noindent (ii) For any $j$ define 
$$Z^j_f=\{ \omega \in \Omega^j: \  \ \dd \omega=0 \text{ and } \dd f \wedge \omega=0 \}$$ and 
$$B^j_f=\{ \omega \in \Omega^j: \  \ \omega = \dd f \wedge \dd (\eta) \text{ for some } \eta \in \Omega ^{j-2} \}.$$
Then it is clear that $B^j_f \subset Z^j_f$ and that there are  natural maps
from the quotient $H^j_f =Z^j_f/B^j_f$  to both cohomology groups $H^j(C^*_f)$ and $H^j(\Omega^*,D_f)$, see also Remark \ref{rksurvive} below.

\end{rk}

\begin{rk}
\label{rkAP}
 Note  that the mixed Hodge structure on 
$H^2(F,\C)_{\ne 1}$ is not pure in general. For a line arrangement, one can use the formulas for the spectrum given in \cite{BS} to study the interplay between monodromy and Hodge structure on 
$H^2(F,\C)_{\ne 1}$.

\end{rk}

\section{Jacobian syzygies and free plane curves}

Consider the graded $S-$submodule $AR(f) \subset S^{3}$ of {\it all relations} involving the derivatives of $f$, namely
$$\rho=(a,b,c) \in AR(f)_q$$
if and only if  $af_x+bf_y+cf_z=0$ and $a,b,c$ are in $S_q$.   The reduced curve $C:f=0$ is a free divisor with exponents $d_1 \leq d_2$ if $AR(f)$ is a free graded $S$-module of rank two, with a basis given by $\rho_1$ and $\rho_2$, syzygies of degree $d_1$ and respectively $d_2$.

To each syzygy $\rho=(a,b,c) \in AR(f)_q$ we associate a differential $2$-form
\begin{equation} 
\label{form}
\omega(\rho)=a \dd y \wedge \dd z -b  \dd x \wedge \dd z+ c \dd x \wedge \dd y\in \Omega^2_{q+2}
\end{equation} 
such that the relation $af_x+bf_y+cf_z=0$ becomes $\dd f \wedge \omega(\rho)=0.$

\begin{ex}
\label{exfree} It is known that the monomial arrangement $\A(m,m,3)$ is free with exponents
$d_1=m+1$ and $d_2=2m-2$, see \cite{Dca}. The forms $\omega(\rho_i)$ for $i=1,2$ corresponding to a basis $\rho_1, \rho_2$ of the $S$-module $AR(f)$ are exactly the forms $\omega_1$ and $\omega _2$ given by
\begin{equation} 
\label{dif1}
\omega_1=(x^{m+1}-2xy^m-2xz^m)\dd y\wedge \dd z -(y^{m+1}-2yz^m-2yx^m)\dd x \wedge \dd z +
\end{equation} 
$$+(z^{m+1}-2zx^m-2zy^m)\dd x\wedge \dd y,$$
and
\begin{equation} 
\label{dif2}
\omega_2=y^{m-1}z^{m-1}\dd y\wedge \dd z -x^{m-1}z^{m-1}\dd x \wedge \dd z+x^{m-1}y^{m-1}\dd x\wedge \dd y.
\end{equation} 
 Note that $\omega_2 \in Z^2_f$, and $\dd \omega_1 \ne 0$ unless $m=3$.
\end{ex}

Hence, up to a shift in degrees, for any polynomial $f$ there is an identification 
$$AR(f)=Syz(f):=\ker \{ \dd f \wedge: \Omega ^2 \to \Omega ^3 \},$$
such that the Koszul relations $KR(f)$ inside $AR(f)$ correspond to the submodule $\dd f \wedge \Omega^1$ in $Syz(f)$. Since $C:f=0$ has only isolated singularities, it follows that $H^1(K^*_f)=0$, i.e.
the following sequence, where the morphisms are the wedge product by $\dd f$, is exact for any $j$
$$ 0 \to \Omega^0_{j-2d} \to \Omega^1_{j-d} \to (\dd f \wedge \Omega^1)_j \to 0.$$
In particular, one has
\begin{equation} 
\label{dimKR1}
\dim (\dd f \wedge \Omega^1)_j =0 \text{  for } j \leq d,
\end{equation} 
and
\begin{equation} 
\label{dimKR2}
\dim (\dd f \wedge \Omega^1)_j =3 {j-d+1 \choose 2} \text{  for } d <j < 2d.
\end{equation} 
On the other hand, we get epimorhisms
\begin{equation} 
\label{epi1}
AR(f)_{q-2}=Syz(f)_q \to H^2(K^*_f)_q,
\end{equation} 
for any $q$. Note that 
$Z^2_f=\ker \{ \dd : Syz(f) \to \Omega ^3 \},$
and there is an obvious morphism
\begin{equation} 
\label{epi2}
H^2_{f,td+k}  \to E_{\infty}^{1-t,t}(f)_k.
\end{equation} 

\begin{rk}
\label{rksurvive} (i) The morphism \eqref{epi2} is not necessarily injective. For the Hessian arrangement considered in subsection \ref{HA}, a direct computation using Singular shows that
$$\dim H^2_{f,16} =2  \text { and }  \dim E_{2}^{0,1}(f)_6=\dim E_{\infty}^{0,1}(f)_6=1.$$
(ii) It is clear that a non-zero  element $[\omega] \in E_1^{1-t,t}(f)_k=H^{2}(K^*_f)_{td+k}$, which has a lifting 
$\omega \in Z^2_f$,
 will survive (i.e. will give rise to a non-zero element) in $E_{\infty}^{1-t,t}(f)_k$. However, two distinct liftings can give rise to the same element in $E_{\infty}^{1-t,t}(f)_k$ in view of (i).
\end{rk}

\begin{rk}
\label{rkdim}
When $C:f=0$ is a free divisor with exponents $d_1 \leq d_2$,
then one clearly has $\dim AR(f)_k=0$ for $k<d_1$, 
\begin{equation} 
\label{dim1}
\dim AR(f)_k={k-d_1+2 \choose 2}
\end{equation} 
for $d_1 \leq k <d_2$, and
\begin{equation} 
\label{dim2}
\dim AR(f)_k={k-d_1+2 \choose 2}+{k-d_2+2 \choose 2}
\end{equation} 
for $k \geq d_2$. For non-free divisors, the computation of the dimensions $\dim AR(f)_k$ (resp. of a basis for $AR(f)_k$) is more complicated. These dimensions give an indication on the size of the linear systems to be solved in the algorithm described in the next section. Sometimes these computation can be avoided, as in the case of Zariski sextic below \ref{exZariski}.

\end{rk}

\section{The algorithm for free curves, nearly free curves, and beyond}
Assume first that $C:f=0$ is a free curve having only weighted homogeneous singularities, with $f$ a reduced, homogeneous polynomial of degree $d$. Let $\alpha_{max}=\frac{q_0}{d}$ be the maximal spectral number of the singularities of $C$ which can be written in this form, i.e. a rational number with denominator $d$. Note that always $q_0<2d$, since all the spectral numbers of plane curve singularities are contained in the interval $(0,2)$.

Assume that $d_1 \leq d_2$ are the exponents of $C$ and let $\rho_1$ and $\rho_2$ be a basis of the $S$-module $AR(f)$ with $\deg \rho_l=d_l$ for $l=1,2$. Then $AR(f)_j=0$ for $j<d_1$ and $AR(f)_j$ for $j\geq d_1$ (resp. $j\geq d_2$) has a basis as a $\C$-vector space obtained by taking all the products $\mu \cdot \rho$ between a monomial $\mu$ in $x,y,z$ and $\rho=\rho_1$ (resp. $\rho =\rho_1$ or $\rho=\rho_2$) having the total degree $j$.
The corresponding differential forms $\mu \cdot \omega(\rho)$ form a basis as a $\C$-vector space for $Syz(f)_q$ where $q=j+2$. Consider now the composition map
\begin{equation} 
\label{Dif1}
\delta_q: Syz(f)_q \xrightarrow{\dd} \Omega^3_q \to H^3(K^*_f)_q
\end{equation} 
where the second map is the canonical projection. In down-to earth terms, we have
\begin{equation} 
\label{Dif1.5}
\delta_q( \omega(\rho)) = \delta_q( a \dd y \wedge \dd z -b  \dd x \wedge \dd z+ c \dd x \wedge \dd y=[a_x+b_y+c_z] \in M(f)_{q-3}.
\end{equation} 
Here $M(f)=S/J_f$ is the Milnor algebra of $f$, with $J_f=(f_x,f_y,f_z)$ the Jacobian ideal of $f$.

Using a computer algebra software as CoCoA
\cite{Co} or Singular \cite{Sing}, one can compute 
\begin{equation} 
\label{Dif2}
\kappa_q := \dim \ker \delta_q,
\end{equation} 
for $d_1+2 \leq q \leq q_0<2d$. Then we compute the difference
\begin{equation} 
\label{Dif3}
\epsilon_q= \kappa _q -\dim (\dd f \wedge \Omega^1)_q
\end{equation} 
using the formulas \eqref{dimKR1} and \eqref{dimKR2}. Since $(\dd f \wedge \Omega^1)_q$ is obviously contained in $\ker \delta_q$, it follows that
\begin{equation} 
\label{Dif4}
\epsilon_q= \dim E^{1-t,t}_2(f)_k
\end{equation} 
where $k \in [1,d]$, $q-k$ divisible by $d$ and $t=(q-k)/d$. Moreover, Theorem \ref{thminj} implies that these are the only terms $E^{1-t',t'}_2(f)_k$ which might be non-zero.

\begin{rk}
\label{rkKR}
When $C$ has only weighted homogeneous singularities, to test if a $2$-form  $\omega(\rho) \in Syz$ as in \eqref{form} belongs to 
$\dd f \wedge \Omega^1$, one uses the description of such forms in terms of the saturation $I_f$ of the Jacobian ideal $J_f$ given in \cite{BJ}, \cite{CoxS}, \cite{DStJSC}. The fact that $C$ is free is equivalent to $I_f=J_f$, and hence in this case one has
\begin{equation} 
\label{KR}
\omega(\rho) \in \dd f \wedge \Omega^1 \text { if and only if } a,b,c \in J_f.
\end{equation}

\end{rk}

\begin{rk}
\label{rkNF0} Assume now that $C:f=0$ is a nearly free curve, with exponents $(d_1,d_2)$.
This means that there is a minimal set of homogeneous generators $\rho_1, \rho_2, \rho_3$ for the $S$-module $AR(f)$, with degrees $d_1$, $d_2$ and $d_2$ respectively, and the module of second order syzygies is spanned by a relation of the form
$$a\rho_1+\ell_2 \rho_2 +\ell_3 \rho_3=0,$$
where $a$ is a homogeneous polynomial of degree $d_2-d_1+1$, and $\ell_2$ and $\ell_3$
are linearly independent linear forms, see \cite{DStNF}. We can assumme that
$$\ell_3(x,y,z)=\alpha x + \beta y +\gamma z \text { with }  \gamma \ne 0.$$
Indeed, by permuting $x,y,z$ if necessary, we can always realize this condition.
With this assumption, one has exactly as above $AR(f)_j=0$ for $j<d_1$ and $AR(f)_j$ for $d_2 > j\geq d_1$  has a basis as a $\C$-vector space obtained by taking all the products $\mu \cdot \rho$ between a monomial $\mu$ in $x,y,z$ and $\rho=\rho_1$ having the total degree $j$.
For $j \geq d_2$, $AR(f)_j$ has a basis as a $\C$-vector space obtained as follows: 
take all the products $\mu \cdot \rho$ between a monomial $\mu$ in $x,y,z$ and $\rho=\rho_1$ or $\rho=\rho_2$ having the total degree $j$, and then add all the products $\nu \rho_3$ with total degree $j$,
where $\nu$ is a monomial only in $x,y$.

With this modification, the above algorithm works exactly as for the free curves, i.e. if all the singularities are weighted homogeneous, then very likely the spectral sequence will degenerate at $E_2$, see Proposition \ref{prop9cusp}. In any case, we can run the above algorithm for $q =j+2 \leq d-1$ and get a precise information on all the roots of the Alexander polynomial $\Delta_C(t)$ by Theorem  \ref{thmEV}, but maybe not on their multiplicities.

\end{rk}

\begin{rk}
\label{rkNF} Assume now that $C:f=0$ is neither a free curve, nor a nearly free one, and that maybe there are non weighted homogeneous singularities on $C$. If among a minimal set of generators for $AR(f)$ provided by Singular or CoCoA, there is only one  syzygy, say $\rho_1$, whose degree  is at most $d-3$, we can run the above algorithm for $q =j+2 \leq d-1$ and get a precise information on all the roots of the Alexander polynomial $\Delta_C(t)$ by Theorem  \ref{thmEV}, but maybe not on their multiplicities. 
On the other hand, if the Alexander polynomial $\Delta_C(t)$ is known, then this approach will lead to explicit bases for the monodromy eigenspaces, see Propositions \ref{propZariski6},  \ref{propZariski12} below.

Similarly, if among a minimal set of generators for $AR(f)$, there are only two  syzygies, say $\rho_1$ and $\rho_2$ whose degrees are at most $d-3$, then one can show that $\rho_1$ and $\rho_2$ can be chosen such that there is no second order syzygy
$g\rho_1+h \rho_2=0$, 
and hence the same remark applies.

\end{rk}

\section{The case of line arrangements}

For a line arrangement $\A$ in $\PP^2$, it  is a major open question whether the monodromy operator $h^1:H^1(F,\C) \to H^1(F,\C)$ is combinatorially determined, i.e. determined by the intersection lattice $L(\A)$, see \cite{OT}. The $1$-eigenspace $H^1(F,\C)_1=H^1(U,\C)$ is known to have dimension $d-1$ and to be  a pure Hodge structure of type $(1,1)$.

 Several interesting examples have been computed by D. Cohen, A. Suciu, A. M\u acinic, S. Papadima, see \cite{CS}, \cite{S1}, \cite{MP}, \cite{S2}.
When the line arrangement $\A$ has only double and triple points, then a complete positive answer is given by S. Papadima and A. Suciu in \cite{PS}. 

\subsection{The monomial arrangement $\A(m,m,3)$}

The case of complex reflexion arrangements is discussed in \cite{MPP}, where the authors prove in particular the following result.

\begin{thm}
\label{thmMPP}
Consider the monomial arrangement 
$$\A(m,m,3): f=(x^m-y^m)(x^m-z^m)(y^m-z^m)=0$$
for $m\geq 2$ and denote by $F(m,m,3)$ the corresponding Milnor fiber.
Let $p$ be a prime number. Then the monodromy operator
\begin{equation} 
\label{mono}
h^*: H^1(F(m,m,3),\C) \to H^1(F(m,m,3),\C)
\end{equation} 
has eigenvalues of order $p^s$ if and only if $p=3$. Moreover, for $p=3$, if we denote the multiplicity of such an eigenvalue by $e^s_3(\A(m,m,3))$, then $e^s_3(\A(m,m,3)) \leq 2$ if $m$ is divisible by $3$, and $e^s_3(\A(m,m,3))\leq 1$
otherwise. Moreover, for $s=1$, both inequalities become equalities.
\end{thm}
For an eigenvalue of order $p^s$, the papers \cite{PS}, \cite{MPP} give also upper bounds for the corresponding multiplicities for any line arrangement, see also \cite{PB}, \cite{BY}. However, for the monomial arrangements, the existence of eigenvalues of order $6$ (resp. $9$) for  $\A(6,6,3)$ (resp. for $\A(9,9,3)$), seems impossible to be decided by the methods used in these papers: the modular inequalities and computations with the associated Orlik-Solomon algebras.

To state  our results, recall the differential $2$-forms associated with the monomial arrangement 
$\A(m,m,3)$ introduced in \eqref{dif1} and \eqref{dif2}.
When $m=3m'$, we also define $\omega'_1=x^{m'-1}y^{m'-1}z^{m'-1}\omega_1$. 
Let  $\iota: F \to \C^3$ denote the inclusion of the Milnor fiber $F=F(m,m,3)$ into $\C^3$ and $\Delta: \Omega^j \to \Omega ^{j-1}$ denote the contraction with the Euler vector field

\begin{thm}
\label{thm1}
Let $\lambda=\exp(2\pi i/3)$ and assume $m\geq 3$. Then the eigenspace 
$$H^{1,0}(F(m,m,3),\C)_{\lambda}=H^1(F(m,m,3),\C)_{\lambda}$$
 is spanned by the differential form $\alpha=\iota^* \Delta (\omega_2)$ when $m$ is not divisible by $3$, and by the forms $\al$ and 
$\be=\iota^* \Delta (\omega' _1)$ when $m$ is divisible by $3$.
\end{thm}
A similar result for $m=2,3$ is proved by Nancy Abdallah in \cite[Section 5.5]{Ab}.
Note that $H^{0,1}(F(m,m,3),\C)_{\overline{ \lambda}}=H^1(F(m,m,3),\C)_{\overline{ \lambda}}$
is just the complex conjugate of $H^1(F(m,m,3),\C)_{\lambda}$, and hence it has a basis given by
the complex conjugate form $\overline{ \alpha}$, resp. forms $\overline{ \alpha}$ and $\overline{ \beta}$.
In principle, it is possible find a basis for $H^{0,1}(F(m,m,3),\C)_{\overline{ \lambda}}$ in terms of polynomial differential forms, but the formulas may become very complicated. One  example where the formulas remain simple is provided by the Zariski sextic with six cusps in Proposition \ref{propZariski6}.

\proof
Assume first that $m$ is not divisible by $3$. Then Theorem \ref{thmMPP} implies that 
$\dim H^1(F(m,m,3),\C)_{\lambda}=1$. The differential form  $\omega _2$ introduced in \eqref{dif2}
belongs to $Syz(f)_{2m}=H^2(K^*_f)_{2m}=E_1^{1,0}(f)_{2m}$, it is clearly non-zero and it has a lifting to $Z^2_f$, given by the form $\omega _2$ itself. By Remark \ref{rksurvive},  $\omega _2$  gives rise to a non-zero element $\alpha$ in $H^1(F(m,m,3),\C)_{\lambda}$ since 
$$\exp(-2\pi i 2m/(3m))=\exp(2\pi i/3)=\lambda.$$
 The formula for this element $\alpha$ comes from the equation \eqref{iso1}.
When $m$ is divisible by $3$, say $m=3m'$, then the differential form  $\omega' _1$ 
belongs also to $Syz(f)_{2m}=H^2(K^*_f)_{2m}=E_1^{1,0}(f)_{2m}$, it is clearly independent of $\omega_2$  in $E_1^{1,0}(f)_{2m}$ since the arrangement is free, and it has a lifting to $Z^2_f$, given by the form $\omega' _1$ itself. 

It remains to explain why the forms $\alpha$ and $\beta$ have Hodge type $(1,0)$. Since the Hodge filtration $F^p$ is contained in the pole order filtration $P^p$, namely $F^p \subset P^p$ for any integer $p$, see \cite[Equation (3.1.3)]{ DS1}, it follows that
$$ H^1(F(m,m,3),\C)_{\lambda}=F^0H^1(F(m,m,3),\C)_{\lambda}=P^0H^1(F(m,m,3),\C)_{\lambda}$$
and
$$F^1H^1(F(m,m,3),\C)_{\lambda} \subset P^1H^1(F(m,m,3),\C)_{\lambda}.$$
This yields an epimorphism
$$Gr_F^0H^1(F(m,m,3),\C)_{\lambda} \to Gr_ P^0H^1(F(m,m,3),\C)_{\lambda}=E^{0,1}_{\infty}(f)_{2m}.$$
On the other hand, we clearly have 
$$\dim E^{0,1}_{\infty}(f)_{2m}=\dim H^1(F(m,m,3),\C)_{\lambda} \geq \dim Gr_F^0H^1(F(m,m,3),\C)_{\lambda}.$$
This implies that the above epimorphism is in fact an isomorphism, and hence 
$F^1H^1(F(m,m,3),\C)_{\lambda} = P^1H^1(F(m,m,3),\C)_{\lambda}.$
Note that 
$$\alpha, \beta \in  E_{\infty}^{1,0}(f)_{2m} =  P^1H^1(F(m,m,3),\C)_{\lambda},$$  which
 completes the proof of Theorem \ref{thm1}.

\endproof

For $m=2,3$ the monomial arrangement is well understood, see \cite[Section 5.5]{Ab}, \cite{S1}, \cite{S2}, so from now on we assume $m\geq 4$.
We have already seen that the line arrangement $\A(m,m,3)$, $m\geq 4$, is free with exponents
$d_1=m+1$ and $d_2=2m-2$. Hence, for $j <d_1$ one has $AR(f)_j=0$, while for
$d_1 \leq j <d_2$ one has
\begin{equation} 
\label{E1}
AR(f)_j=S_{j-d_1}\rho_1=\{(ga_1,gb_1,gc_1) \ \ : \ \ g \in S_{j-d_1}\},
\end{equation}
where one has $a_1=x^{m+1}-2xy^m-2xz^m$, $b_1=y^{m+1}-2yz^m-2yx^m$ and $c_1=
z^{m+1}-2zx^m-2zy^m.$
It follows that in these cases the morphism $\delta _q$ from \eqref{Dif1}, \eqref{Dif1.5}, for $q=j+2$, can be identified with the morphism $\phi_j: S_{j-d_1} \to M(f)_{j-1}$ given by
\begin{equation} 
\label{E2}
g \mapsto [ (ga_1)_x+(gb_1)_y+(gc_1)_z]. 
\end{equation}
It is easy to prove the following result.

\begin{lem}
\label{lemE1}
The morphism $\phi_j$ is injective for any $j$ with $d_1 \leq j <d_2$.
\end{lem}

\proof
For $g \in S_{j-d_1}$ we have
$$\phi_j(g)=x^m(3xg_x+rg)+y^m(3yg_y+rg)+z^m(3zg_z+rg),$$
with $r=m-3-2(j-d_1)$. Note that $j-d_1<d_2-d_1=m-3<m$ and $M(f)_{j-1}=S_{j-1}$.
Since $x^m,y^m,z^m$ is a regular sequence in $S$, the equality $\phi_j(g)=0$ implies
$$3xg_x+rg=3yg_y+rg=3zg_z+rg=0.$$
This yields $(j-d_1+r)g=0$. Note that $j-d_1+r=m-3-(j-d_1)=d_2-j >0$, and hence $g=0$.

\endproof

It follows that it is enough to consider the case $j \geq d_2$. Then
\begin{equation} 
\label{E3}
AR(f)_j=S_{j-d_1}\rho_1\oplus S_{j-d_2}\rho_2=\{(ga_1+ha_2,gb_1+hb_2,gc_1+hc_2)   :  g \in S_{j-d_1},  \ h \in S_{j-d_2}\},
\end{equation}
where $a_2=y^{m-1}z^{m-1}$, $b_2=x^{m-1}z^{m-1}$, $c_2=x^{m-1}y^{m-1}$. In these cases the morphism $\delta _q$ from \eqref{Dif1}, \eqref{Dif1.5}, for $q=j+2$, can be identified with the morphism $$\phi_j: S_{j-d_1} \oplus S_{j-d_2}\to M(f)_{j-1}$$ given by
\begin{equation} 
\label{E4}
(g,h) \mapsto [ (ga_1+ha_2)_x+(gb_1+hb_2)_y+(gc_1+hc_2)_z]. 
\end{equation}

\begin{rk}
\label{rklimit}
In view of Example \ref{exspec}, it is enough to consider $j \in [d_2, 6m-6]$ if we want to check the $E_2$-degeneracy of the spectral sequence. On the other hand, if we want just to determine all the eigenvalues, in view of Theorem \ref{thmEV}, it is enough to consider $j \in [d_2, 3m-3]$.
\end{rk}

We give below the values of the differences $\epsilon_q$, for $q=j+2$, as defined in \eqref{Dif3}, and obtained using Singular.
\begin{thm}
\label{thmE1}
Assume that $4 \leq m \leq 25$. Then $\epsilon _q=0$, except for 
\begin{enumerate}

\item $q=2m$ and $q=4m$, when
$\epsilon _q=e^1_3(\A(m,m,3))$, i.e. $1$ if $m$ is not divisible by $3$ and $2$ otherwise. 

\item $q=3m$,  when $\epsilon _q=3m-1.$

\end{enumerate}

In particular one has the following, 

\begin{enumerate}
\item[(i)]  Conjecture \ref{conj} holds for the line arrangement $\A(m,m,3)$ for $m \in [4,25]$.

\item[(ii)]  For $m \in [4,25]$, the monodromy operator \eqref{mono} has only cubic roots of unity as eigenvalues.

\end{enumerate}

\end{thm}

If the arrangement  $\A(m,m,3)$ were deformable into a real arrangement, the claim (ii) would follow from \cite{Yo0}, but this does not seem to be the case. The arrangement $\A(m,m,3)$ is definitely not deformable to a real line arrangement, as can be seen by an argument similar to \cite[Proposition 4.6]{NY}.

\begin{ex}
\label{exaction2}
Let $G$ be a cyclic group of order $2m$, whose generator $T$ acts on $\C^3$ by
$$T(x,y,z)=(\xi y, \xi x, \xi z),$$
where $\xi$ is a primitive root of unity of order $2m$. Then the Milnor fiber $F(\A(m,m,3))$ is clearly $G$-invariant, and hence there is an induced $G$-action on $F(\A(m,m,3))$.
Using Theorem \ref{thm1}, it follows that
$T^*=- Id$ on $H^1(F(\A(m,m,3)),\C)_{\ne 1}$.
\end{ex}

\subsection{The Hessian arrangement } \label{HA}

The Hessian arrangement is given by the equation 
$$\A: f=xyz\left( (x^3+y^3+z^3)^3-27x^3y^3z^3\right )=0,$$
and consists of all the four  singular members of the pencil $\PPP : uc_1+vc_2$, where $$c_1=x^3+y^3+z^3 \text { and }c_2=xyz.$$
Each of these four members is a triangle.
This arrangement plays a key role in the theory of line arrangements, as it is the only known $4$-net, see \cite[Example 2.15]{S2},  \cite[Example 8.7]{PS}. Its monodromy is computed in \cite[Remark 3.3 (iii)]{BDS},
\cite[Theorem 1.7]{PS}

Here we describe bases for the nontrivial eigenspaces $H^1(F,\C)_{\lambda}$, with
$\lambda \ne 1$. To state  our results, let us introduce some differential $2$-forms associated with the Hessian arrangement. Define
\begin{equation} 
\label{difHA1}
\omega_1=-x(y^3-z^3)\dd y\wedge \dd z -y(x^3-z^3 )\dd x \wedge \dd z +z(y^3-x^3)\dd x\wedge \dd y=dc_1 \wedge dc_2,
\end{equation} 
and
\begin{equation} 
\label{difHA2}
\omega_2=A\dd y\wedge \dd z -B\dd x \wedge \dd z +C\dd x\wedge \dd y,
\end{equation}
where
$$A=x^7+3x^4y^3+3xy^6+13x^4z^3-81xy^3z^3+34xz^6 ,$$
$$B=y^7+16y^4z^3-44yz^6,$$
and
$$C= -10x^6z-21x^3y^3z-13y^6z-31x^3z^4+47y^3z^4+z^7)\dd x\wedge \dd y.$$ 
It is known that $\A$ is a free arrangement with exponents $(4,7)$, see \cite{JV}, \cite{Dca}, and the generating syzygies correspond to the forms $\omega_1$ and $\omega_2$.

\begin{thm}
\label{thmHA}

For the  Hessian arrangement, one has the following. 
\begin{enumerate}

\item The eigenspace $H^{1,0}(F,\C)_{-1}$ is spanned by $\alpha=\iota^* \Delta (\omega_1)$.

\item The eigenspace $H^{1,0}(F,\C)_{i}=H^{1}(F,\C)_{i}$ is spanned by $\beta_1=\iota^* \Delta (c_1\omega_1)$ and
 $\beta _2=\iota^* \Delta (c_2 \cdot \omega_1)$.
\end{enumerate}

\end{thm}

\proof Exactly as the proof of Theorem \ref{thm1}.

\endproof

\begin{ex}
\label{exaction1}
Let $G=\Sigma_3$, the permutation group acting on $\C^3$ by permuting the coordinates $x,y,z$. Then it is clear that the Milnor fiber $F$ is $G$-invariant, hence we get an induced action of $G$ on $F$. Since the differential form $\omega_1$ is clearly invariant under $G$ (it is enough to check for a transposition $x \mapsto y,$ $y \mapsto x$, $z \mapsto z$), and the cubic polynomials $c_1$ and $c_2$ are symmetric, it follows that the action of $G$ on 
$H^1(F,\C)_{\ne 1}$ is trivial. Note that the $G$-action on $H^1(F,\C)_{1}$ is not trivial, and can be easily described.
\end{ex}

\section{Further examples}

Recall that for a plane curve $C:f=0$ of degree $d$, the integer $k$ takes values in the interval $[1,d]$.

\subsection{Free curves with non weighted homogeneous singularities} \label{freecurve}

 First we consider the rational cuspidal curve
$$
C:f=(xz-y^2)^3-x^2 y^4=0,
$$
which is  free  with exponents $d_1=2$, $d_2=3$.
It has a weighted homogeneous $E_6$-singularity located at the point $[1:0:0]$,  and a non-weighted homogeneous  $E_{14}$-singularity located at  $[0:0:1]$.  This latter singularity has the same topological type as the singularity
$u^3-v^8$, it has  Milnor number~$\mu(E_{14})=14$, and its Tjurina number is~$\tau(E_{14})=13$. It follows from the discussion in \cite{AD} that the fundamental group of the complement  coincides
with $\Z/2*\Z/3$ and the Alexander polynomial  of $C$ is given by
$$\Delta_C(t):=\det (t\cdot Id -h^*|H^1(F,\C))= t^2-t+1.$$
By applying our computational approach to this curve, we obtain the next result.

\begin{prop}
\label{propNWH} 
 For the free curve
$
C:f=(xz-y^2)^3-x^2 y^4=0,
$ the following hold.

\begin{enumerate}

\item The generating syzygies $\rho_1$ and $\rho_2$ for $AR(f)$ are given by
$$\rho_1=(2x^2, -xy,-2y^2-2xz) \text{ and } \rho_2=(6y^3-18xyz,3y^2z+3xz^2,4y^3+24yz^2).$$

\item $\dd (y\omega(\rho_1))=0.$ In particular $\alpha=\iota^*(\Delta(y\omega(\rho_1))$ spans the 1-dimensional vector space 
$H^{1,0}(F,\C)_{\lambda}=H^1(F,\C)_{\lambda},$
 for $\lambda=\exp(\pi i/3).$

\item   $\dim E_2^{1,0}(f)_k=1$ for $k=5$ and vanishes otherwise. 

\item   $\dim E_2^{0,1}(f)_k=1$  for $k=1, 4,5,6$ and and vanishes otherwise.

\item   $\dim E_2^{1-t,t}(f)_k=1$ for any $t \geq 2$ and any $k \in [1,6].$

\end{enumerate}

\end{prop}
Notice that the spectral sequence $E_*(f)$ does not degenerate in this case, see \cite[Theorem 5.2]{DS1}, and the dimensions 
$\epsilon _q$ from \eqref{Dif3} verify
$$\epsilon _q=1=\mu(E_{14})-\tau(E_{14}),$$
for $q$ large, in fact $q \leq 10$, as predicted by the theory, see \cite[Proposition 3.7 and Theorem 3.9]{Dcomp}.
On the other hand, the claim (3) implies the formula for $\Delta_C(t)$ given above, in view of Theorem \ref{thmEV}.

\medskip

\noindent {\bf Proof of Corollary \ref{corAP}.}
Now we explain along the same lines as above the proof of Corollary \ref{corAP}. 
For any integer $m\geq 1$, the generating syzygies  can be chosen as

\noindent $\rho_1=(0,2y^{m+1}z^{m-1}, x^{2m}-3y^mz^m)$
and $\rho_2=(A,B,C ),$
where the entries $A,B,C$ are given by 
$A=2x^{2m}y^{m-1}-2y^{2m-1}z^m,$
$B=10x^{3m-1}-8x^{2m-1}y^m+30x^{m-1}y^mz^m$ and $C=8x^{2m-1}y^{m-1}z-45x^{m-1}y^{m-1}z^{m+1}.$

When $m$ is odd,
$1 \leq m \leq 19$, a direct computation using Singular shows that
$ E_2^{1,0}(f)_k=0$ for any $k \in [1,d]$, with $d=5m$, which gives the claim in this case via Theorem \ref{thmEV}. When $m$ is even, $2 \leq m \leq 20$, the same computation gives
$ E_2^{1,0}(f)_k=0$ for any $k \in [1, d/2]$ and $\dim E_2^{1,0}(f)_k=1$ for any $k \in (d/2, d]$,
which allow us to conclude again via Theorem \ref{thmEV}. Note also that the curve $C_{5m}$ has two singular points, namely at $[0:1:0]$ and $[0:0:1]$, and a local analysis at these points using general principles shows that
$$\tau(C_{5m})=19m^2-8m+1 < 21m^2-13m+2 \leq \mu(C_{5m}),$$
for $m \geq 3$. A  computation by Singular in the case $m=2$ gives  $\tau(C_{10})< \mu(C_{10})$.
 In other words, for $m \geq 2$,  the curve  $C_{5m}$ has at least a non weighted homogeneous singularity, and hence the spectral sequence $E_*(f)$ does not degenerate at $E_2$.

\subsection{A nearly free curve: the sextic with 9 cusps} \label{ex9cusps}

Consider the sextic curve
$$C:f=x^6+y^6+z^6-2(x^3y^3+x^3z^3+y^3z^3)=0,$$
having $9$ cusps  and Alexander polynomial $\Delta_C(t)= (t^2-t+1^3,$ see for instance \cite[Corollary 12]{Oka1}. This curve is nearly free with exponents $(3,3)$. More precisely, the generating syzygies are
$$ \rho_1=(y^3-z^3, x^2y,-x^2z),$$
$$ \rho_2=(xy^2, x^3-z^3,-y^2z),$$
$$ \rho_3=(xz^2,-yz^2, x^3-y^3).$$
The generating second order syzygy is $x\rho_1-y\rho_2+z \rho_3=0$. The algorithm described in Remark \ref{rkNF0} gives us the following.
\begin{prop}
\label{prop9cusp}  Let  $C:f=x^6+y^6+z^6-2(x^3y^3+x^3z^3+y^3z^3)=0$ be the above sextic curve with $9$ cusps. Then the following hold.
\begin{enumerate}

\item Conjecture \ref{conj} holds for the sextic curve $C$.

\item   $\dim E_2^{1,0}(f)_k=3$ for $k=5$, $\dim E_2^{0,1}(f)_k=3$ for $k=1$, 
 and all the other terms $E_2^{1-t,t}(f)_k$ vanish.

\item  The $1$-forms $\alpha_j=\iota^*(\Delta(\omega(\rho_j))$ for $j=1,2,3$ span the $3$-dimensional vector space 
$H^{1,0}(F,\C)_{\lambda}=H^1(F,\C)_{\lambda}$, where $\lambda=\exp(\pi i/3).$

\end{enumerate}
\end{prop}

\subsection{A Zariski sextic with six cusps} \label{exZariski}

Consider the sextic curve
$$C:f=(x^2+y^2)^3+(y^3+z^3)^2=0,$$
having $6$ cusps on a conic and Alexander polynomial $\Delta_C(t)= t^2-t+1,$ see for instance \cite[Theorem 6.4.9]{D1}. This curve is neither free, nor nearly free: the module $AR(f)$ has a minimal set of generators consisting of one syzygy of degree $3$, namely
$$\rho_1= (yz^2, -xz^2, xy^2),$$
and $3$ other syzygies of degree $5$, among which
$$\rho_2= (y^3z^2 + z^5, 0, -x^5 - 2x^3y^2 - xy^4).$$
 It is clear that $\dd (\omega(\rho_1)) =\dd (\omega(\rho_2))       =0$, hence we get the following, using the algorithm described in Remark \ref{rkNF}.

\begin{prop}
\label{propZariski6}  Let $C:f=(x^2+y^2)^3+(y^3+z^3)^2=0$ be the above sextic curve and 
set $\lambda=\exp(\pi i/3).$ Then the following hold.
\begin{enumerate}

\item   $\dim E_2^{1,0}(f)_k=1$ for $k=5$ and vanishes otherwise. 

\item  $\alpha=\iota^*(\Delta(\omega(\rho_1))$ spans the 1-dimensional vector space 
$H^{1,0}(F,\C)_{\lambda}=H^1(F,\C)_{\lambda}.$

\item $\beta=\iota^*(\Delta(\omega(\rho_2))$ spans the 1-dimensional vector space 
$H^{0,1}(F,\C)_{\overline\lambda}=H^1(F,\C)_{\overline \lambda}.$

\end{enumerate}
\end{prop}

Note that  the claim (1) implies the formula for $\Delta_C(t)$ given above, in view of Theorem \ref{thmEV}.

\subsection{A $(3,4)$-torus type curve } \label{torus}

Consider the degree $12$ curve
$$C:f=(x^3+y^3)^4+(y^4+z^4)^3=0,$$
which is  a $(3,4)$-torus type curve having $12$ singularities $E_6$ and Alexander polynomial given by $\Delta_C(t)= \Phi_6(t) \Phi_{12}(t),$ with $\Phi_m(t)$ being the $m$-th cyclotomic polynomial, see for instance \cite[Theorem 31]{OkaS}. This curve is neither free, nor nearly free: the module $AR(f)$ has a minimal set of generators consisting of one syzygy of degree $3$, namely
$$\rho= (y^2z^3, -x^2z^3, x^2y^3),$$
and the other generating syzygies of degree at least $11$.  It is clear that 
$$\dd (\omega) =\dd (\eta_1)=\dd (\eta_2)=0,$$
where $\omega=\omega(\rho)$, $\eta=\omega( (x^3+y^3) \rho)$ and $\omega'=\omega( (y^4+z^4) \rho)$.
Set $\lambda=\exp(\pi i/6)$ and note that $\lambda^2$ and $\overline \lambda^2$ are the roots of $\Phi_6(t)$, and $\lambda$, $\lambda^5$, $\overline \lambda$, $\overline \lambda^5$ are the roots of $\Phi_{12}(t)$. The algorithm described in Remark \ref{rkNF} yields the following.

\begin{prop}
\label{propZariski12}  Let $C:f=(x^3+y^3)^4+(y^4+z^4)^3=0$ be the above $(3,4)$-torus type curve. 
 Then the following hold.
\begin{enumerate}

\item   $\dim E_2^{1,0}(f)_k=1$ for $k=7,10,11$ and vanishes otherwise.

\item  $\alpha=\iota^*(\Delta(\omega)$ spans the 1-dimensional vector space 
$H^{1,0}(F,\C)_{\lambda^5}=H^1(F,\C)_{\lambda^5}.$

\item $\alpha'=\iota^*(\Delta(\omega')$ spans the 1-dimensional vector space 
$H^{1,0}(F,\C)_{\overline\lambda}=H^1(F,\C)_{\overline \lambda}.$

\item $\beta=\iota^*(\Delta(\eta)$ spans the 1-dimensional vector space 
$H^{1,0}(F,\C)_{\lambda^2}=H^1(F,\C)_{\lambda^2}.$

\end{enumerate}
\end{prop}

Again  the claim (1) implies the formula for $\Delta_C(t)$ given above, in view of Theorem \ref{thmEV}.

\end{document}